\documentclass[11pt]{article}
\usepackage{authblk}
\usepackage[T1]{fontenc}
\usepackage[utf8]{inputenc}
\usepackage[english]{babel}
\usepackage{lmodern}
\usepackage{mathtools} 
\usepackage[a4paper,DIVcalc]{typearea} 
\typearea[10mm]{16} 

\providecommand{\abstractt}[1]
{{
  \small
  \textbf{Abstract -} #1
}}

\providecommand{\keywords}[1]
{{
  \small 
  \textbf{Keywords -} #1
}}

\providecommand{\AMSClassification}[1]
{{
  \small
  \textbf{AMS Subject classifications -} #1
}}

 \usepackage{comment}
\usepackage{cases}
\usepackage{setspace}
\usepackage{amsthm,amssymb,amsmath,eqnarray}
\usepackage{dsfont}
\usepackage{latexsym}
\usepackage[pdftex]{graphicx}
\usepackage{mathdots}
\usepackage[hyperindex=true, pdftex=true, linkcolor=black,colorlinks=true,citecolor=blue,urlcolor=blue]{hyperref}
\usepackage{amssymb,amsfonts,textcomp}
\usepackage{multirow}
\usepackage{longtable}
\usepackage{xcolor}
\definecolor{LightGray}{gray}{0.9}

\usepackage{pgf, tikz}
\usepackage{authblk}

\title{Symmetry Analysis of Semi-Linear Partial Differential Equations and Forward Backward Stochastic Differential Equations  }
\author{OUKNINE Anas, LESCOT Paul}

\affil{Laboratoire de Math\'ematiques Rapha\"el Salem, University of Rouen, UMR CNRS 6085, Avenue de l'Universit\'e, 76801 Saint Etienne du Rouvray, France\\
Email: \texttt{anas.ouknine@univ-lemans.fr}\\
Email: \texttt{paul.lescot@univ-rouen.fr}}
\date{\today}

\newtheorem{thm}{Theorem}[section]
\newtheorem{defin}[thm]{Definition}
\newtheorem{prop}[thm]{Proposition}

\newtheorem{remark}[thm]{Remark}
\newtheorem{corollary}[thm]{Corollary}

\numberwithin{equation}{section}
\theoremstyle{remark}

\usepackage{cases}
\usepackage{empheq}
\usepackage{hyperref}

 \allowdisplaybreaks[1]

\begin{document}

\maketitle


\noindent \abstractt{We examine the Lie symmetries of a semi-linear partial differential equations and their connections to the analogous symmetries of the forward-backward stochastic differential equations (FBSDEs), established through the generalized Feynman-Kac formula.}
\newline

\noindent \keywords{Forward Backward Stochastic Differential Equations, Semi-linear Partial Differential Equation, Symmetry, Lie Algebra, Vector field}\newline

\noindent\AMSClassification{34A26, 91B28 ,60H10, 60H30, 58D19}

\section{Introduction}
Symmetry analysis of differential equations is a powerful tool in the study of ordinary differential equations (ODEs) and partial differential equations (PDEs). This theory is a very general and useful tool for finding analytical solutions of large classes of differential equations. A comprehensive overview of this method can be found in \cite{OLV86,ovsiannikov2014group,ibragimov1995crc,blu}. It is thus natural to think they could be useful also in the study of stochastic differential equations (SDEs). This observation is certainly not new, and there is now a substantial body of literature dedicated to the study of symmetries of SDEs \cite{gaeta2017random,Giuseppe-Gaeta_2000,kozlov2010symmetries,kozlov2020symmetries,unal2003symmetries,srihirun2007definition,gaeta2019w,Giuseppe-Gaeta_1999}.\\
Simply speaking, the symmetry method involves a group of transformations (known as Lie group), that maps a solution of a given system of equations to another solution of the same system. Since then many papers contributed to this area.\\
In 1999, Gaeta applied the symmetry method to SDEs, involving transformations of the Brownian motion \cite{Giuseppe-Gaeta_1999}, as well as to the associated Kolmogorov Forward Equations (Fokker-Planck equations). He also analyzed the relationship between the symmetries of the SDEs and the Fokker-Planck equations, identifying conditions under which the symmetry of one implies the symmetry of the other.\\
A very recent application of the symmetry method to backward stochastic differential equations (BSDEs) and forward-backward stochastic differential equations (FBSDEs) was conducted by Na Zhang and Guangyan Jia \cite{sym11091153,ZHANG2021105527}, which will be the focus of our interest in this article. \\
FBSDEs hold particular importance due to their profound connections with PDEs, especially semi-linear PDEs of type:
\begin{align}
    \begin{cases}
        &\partial_t u(t, x)+[\mathcal{L}_{t} u](t, x)+g\left(t, x, u(t, x), \partial_x u(t, x) \sigma(t, x)\right)=0,(t, x) \in[0, T) \times \mathbb{R}, \\
        &u(T, x)=H(x), \quad x \in \mathbb{R}. \label{SemiLinear-PDE}\\
    \end{cases}
\end{align}
 \begin{equation}
     \text{where}\quad\mathcal{L}u=b(t,x)\partial_x u(t, x)+\frac{1}{2}\sigma^2(t,x)\partial^{2}_{xx} u(t, x)\label{inf gen of sde}
 \end{equation}
is the infinitesimal generator of the SDE 
\begin{equation}
    dX_t = b(t, X_t) dt + \sigma(t, X_t) dW_t.\label{SDE**}
\end{equation}
 This relationship is famously illuminated by the Feynman-Kac formula, which acts as a bridge between stochastic processes and PDEs by providing a probabilistic representation of the solution to certain PDEs. It achieves this by expressing the solution as the expectation of functional of stochastic processes. This foundational principle naturally extends to FBSDEs and semi-linear PDEs.
Specifically, the solutions of FBSDEs can be shown to correspond to the solutions of the associated semi-linear PDEs. Understanding this connection provides significant enrichment to both the theory and practical applications of stochastic processes and PDEs. It allows us to utilize the tools and techniques from one domain to solve problems in the other, creating a powerful framework for addressing complex mathematical and real-world challenges. This interplay not only enriches our theoretical insights but also enhances our ability to model and solve various phenomena in fields such as finance, physics, and engineering.\\
Motivated by Gaeta's and Kozlov's results \cite{Giuseppe-Gaeta_1999,kozlov2021symmetries}, we aim to explore the symmetries of FBSDEs and semi-linear PDEs, and investigate the relation between them. 
This paper is organized as follows: Section 2 presents the fundamental concepts of FBSDEs, uncluding the  conditions for existence and uniqueness of solutions, as well as their connection to the semilinear PDE through the generalized (non linear) Feynman-Kac formula. Additionally, we introduce the definition of of symmetry of FBSDEs, following the approach of Gaeta and Na Zhang \cite{Giuseppe-Gaeta_1999,sym11091153}.
Section 3 focuses on computing the symmetries of semilinear PDEs using Olver's prolongation method \cite{OLV86}. In Section 4, we compare the symmetries of FBSDEs with those of semilinear PDEs under certain restrictions on the generator of the FBSDE. Finally, Section 5 provides an example to illustrate the results presented in Section 4.
\section{Some Facts about FBSDE and their Symmetries}
\subsection{FBSDE}
Fix a terminal time $T>0$, and let $B_t$ a standard Brownian motion on $\mathbb{R}$ defined on a probability space $(\Omega,(\mathcal{F}_t)_{0\leq t\leq T},\mathcal{F},\mathbb{P})$, and $(\mathcal{F}_t)_{0\leq t\leq T}$ is the filtration generated by $B_t$. Consider an uncoupled FBSDE, which typically takes the form:
\begin{align}
    \begin{cases}
        dX_t &= b(t, X_t) dt + \sigma(t, X_t) dW_t, \\
        dY_t &= -g(t, X_t, Y_t, Z_t) dt + Z_t dW_t, \label{FBSDE}\\
        Y_T &= H(X_T).
    \end{cases}
\end{align}

\noindent The generator $g$ is defined on $[0,T]\times \mathbb{R}\times \mathbb{R}\times \mathbb{R}$ and is assumed to always be measurable with respect to $\mathcal{B}([0,T])\times \mathcal{B}(\mathbb{R})\times \mathcal{B}(\mathbb{R})\times \mathcal{B}(\mathbb{R})$. Additionally, for any $X$, $Y$ and $Z$ the process $t\longrightarrow g(t,X_t,Y_t,Z_t)$ is progressively measurable. 
The terminal condition $Y_T$ is a $\mathcal{F}_T$-measurable random variable, typically denoted by $\xi$. The solution $(X_t,Y_t,Z_t)_{t\in [0,T]}$ consists of three stochastic processes $(X_t)_{t\in [0,T]}$, $(Y_t)_{t\in [0,T]}$ and $(Z_t)_{t\in [0,T]}$ which are adapted to the filtration $(\mathcal{F}_t)_{t\in [0,T]}$.
It is called uncoupled because the solution
$(Y_t,Z_t)$ of the BSDE in \eqref{FBSDE} does not interfere with the dynamics of the forward SDE in \eqref{FBSDE}. This allows the forward equation to be solved first, and its evolution (solution) can be inserted into the backward equation, which can then be solved as a BSDE. Unlike standard SDEs, BSDE are solved backwards in time. BSDEs naturally arise in various applications such as stochastic control, mathematical finance, and nonlinear Feynman-Kac formulae. BSDEs were introduced by Jean-Michel Bismut in 1973 in the linear case \cite{bismut1973conjugate} and by Étienne Pardoux and Shige Peng in 1990 in the nonlinear case \cite{PARDOUX199055}. The conditions for the existence and uniqueness of BSDEs are varied. For instance, these conditions can include:
\begin{enumerate}
    \item $g$ is Lipschitz continuous in $y$ and $z$, i.e., 
    see for example \cite{PARDOUX199055}.
    \item 
    $g$ is locally Lipschitz in both variables $y$ and $z$, see \cite{10.1214/ECP.v7-1058}.
    \item $g$ of quadratic growth in $z$, for example see \cite{briand2008quadratic,10.1214/12-AOP743,10.1214/aop/1019160253}
    \item $g$ is uniformly continuous and of at most linear growth in $y$, $z$, see \cite{jia2009some}.
\end{enumerate}

\noindent Various assumptions about g may lead to different conditions for $Y_T$, $Y$, and $Z$. Since the symmetry
method only requires the existence of the solutions of BSDEs, we will not specify the conditions. Note that $g$ is usually assumed to be dependent on $\omega \in \Omega$ in the literature on BSDE. However, in this
article, similarly to most applications of symmetry methods in SDEs, we assume $g$ to be non-random (deterministic).

\subsection{Relation between FBSDE and semi-linear PDE}
Consider the following PDE:
\begin{align}
    \begin{cases}
        &\partial_t u(t, x)+[\mathcal{L}_{t} u](t, x)-f(t,x)u(t,x)=0,(t, x) \in[0, T) \times \mathbb{R}, \\
        &u(T, x)=H(x), \quad x \in \mathbb{R}. \label{Feynman Kac-PDE}\\
    \end{cases}
\end{align}
where $\mathcal{L}_t$ is the infinitesimal generator of the SDE \eqref{SDE**} defined in \eqref{inf gen of sde}.
Using Feynman-Kac formula \cite{oksendal2013stochastic}, the solution of the PDE \eqref{Feynman Kac-PDE} is of the form 
\begin{equation}
    u(t,x)=\mathbb{E}_x\left(e^{-\int_t^T f(s,X_s)ds}H(X_T) \right)
\end{equation}
Here, $X_t$ is a solution of the diffusion \eqref{SDE**}. Thus, we have found a probabilistic representation of the solution to the PDE \eqref{Feynman Kac-PDE}. However, if the potential term in the PDE \eqref{Feynman Kac-PDE} is not linear in $u$, and more precisely, if the potential term is of the form  $f(t,x,u,\sigma(t,x)u(t,x))$, we are led to consider a FBSDE as the probabilistic representation of the solution to the PDE with this nonlinear potential. We will further elaborate on this relationship in the following proposition.

\begin{prop}\cite[p.~42]{carmona2016lectures}\label{PDE-FBSDE}
 Let us assume that the function $(t, x) \longrightarrow u(t, x)$ is jointly continuous on $[0, T] \times \mathbb{R}^d$, continuously differentiable in $t$, and twice continuously differentiable in $x$, the first derivative $\partial_x u(t, x)$ being of polynomial growth in $x$. Let us also assume that this function $u$ is a classical solution of the PDE
$$
\left\{\begin{array}{l}
\partial_t u(t, x)+[\mathcal{L} u](t, x)+g\left(t, x, u(t, x), \partial_x u(t, x) \sigma(t, x)\right)=0,(t, x) \in[0, T) \times \mathbb{R}^d, \\
u(T, x)=H(x), \quad x \in \mathbb{R}^d
\end{array}\right.
$$
Then $(\mathbf{Y}, \mathbf{Z})=\left(Y_t, Z_t\right)_{t \in[0, T]}$, defined by
$$
Y_t=u\left(t, X_t\right), \quad Z_t=\partial_x u\left(t, X_t\right) \sigma\left(t, X_t\right), \quad 0 \leq t \leq T
$$
is the unique solution of the BSDE
$$
d Y_t=-g\left(t, X_t, Y_t, Z_t\right) d t+Z_t d W_t, \quad 0 \leq t \leq T, \quad Y_T=H\left(X_T\right),
$$
with respect to the filtration generated by the standard Brownian motion $B_t$.
Here, $X_t$ is the diffusion \eqref{SDE**}
\end{prop}

\subsection{Symmetries of FBSDE}
Let us recall some facts about Lie symmetries of a differential equations. When co
\begin{defin}
    We call projectable symmetries, vector fields of the form
    \begin{equation}
    v_{_{FBSDE}}=\tau(t) \partial_t+\xi(t, x, y) \partial_x+\eta(t, x, y) \partial_u.\label{Projectable VF}
\end{equation}
means transformations of time $\tau$ depends only on $t$.
\end{defin}
\noindent Recall the FBSDE \eqref{FBSDE}. We introduce a projectable vector field \eqref{Projectable VF}
\noindent where $\tau:\mathbb{R}_+\longrightarrow \mathbb{R},\quad \xi:\mathbb{R}_+\times \mathbb{R}\times \mathbb{R}\longrightarrow \mathbb{R}$ and $\eta:\mathbb{R}_+\times \mathbb{R}\times \mathbb{R}\longrightarrow \mathbb{R}$ are smooth. Every vector field generates a one-parameter group of transformations (flow) denoted by $$\psi_\varepsilon(t,x,y):=\exp(\varepsilon v)(t,x,y)=(\Tilde{t},\Tilde{x},\Tilde{y})$$
\begin{equation}
      \text{such that }\quad \frac{d\tilde{x}}{d\varepsilon}=\xi(\tilde{x},\tilde{t},\tilde{y})\,\,,\,\, \frac{d\tilde{t}}{d\varepsilon}=\tau(\tilde{t})\,\,,\,\,\frac{d\tilde{u}}{d\varepsilon}=\eta(\tilde{x},\tilde{t},\tilde{y})\label{Lie's equation(exponentiation)}.
\end{equation}
 
\noindent with initial conditions 
$\tilde{x}(0)=x\,,\,\tilde{t}(0)=t\,,\,\tilde{u}(0)=u.$ For further information on these notations, readers can refer to \cite{OLV86}.

We now provide a definition of the symmetry of FBSDE.
\begin{defin}\label{def of sym of FBSDE}
The vector field \eqref{Projectable VF} is a symmetry of the FBSDE \eqref{FBSDE} if the transformed process $(\tilde{X}_{\tilde{t}},\tilde{Y}_{\tilde{t}},\tilde{Z}_{\tilde{t}})$ is a solution of \eqref{FBSDE}.
    \end{defin}
\begin{remark}
\begin{itemize}
    \item 
    In \eqref{Projectable VF}, we did not include $z$ in the transformations $\xi$ and $\eta$ because $z$ is a process related to $y$, we can obtain the transformation for $z$ simultaneously as we derive the transformation for $y$.
    \item 
    According to the symmetry \eqref{Projectable VF}, if $\tau \neq 0$, this action, hence, transform the Brownian motion $W_t$  to a new Brownian motion called $\tilde{W}_{\tilde{t}}$ defined by:
    \begin{equation}
        \tilde{w}_s=\sqrt{\frac{d \tilde{t}}{dt}} w_s=\sqrt{1+\varepsilon \tau^{'}} w_s \label{new brownian oks}
    \end{equation} this result is known as 'the random time change formula' due to Oksendal in \cite{oksendal2013stochastic,oksendal1985stochastic}. Many authors used this relation to compute symmetry of SDE using different approaches \cite{Giuseppe-Gaeta_1999,Giuseppe-Gaeta_2000,gaeta2019w, gaeta2017random, kozlov2010symmetries,kozlov2020symmetries, unal2003symmetries,srihirun2007definition}.
\end{itemize}    
\end{remark}
Generally for variables $(t,x,y,z)$ and the transformation $\psi_{\varepsilon}$, we denote simply the transformed variables by $(\tilde{t},\,\tilde{x},\,\tilde{y}, \tilde{z})$ or $(\tilde{t}_{\varepsilon},\,\tilde{x}_{\varepsilon},\,\tilde{y}_{\varepsilon}, \tilde{z}_{\varepsilon})$. Let $(x,y,z)$ be a solution of \eqref{FBSDE}. Then definition \ref{def of sym of FBSDE} means that the vector field $v_{_{FBSDE}}$ is a symmetry of \eqref{FBSDE} if and only if the transformed process $(\tilde{x}_{\tilde{t}},\,\tilde{y}_{\tilde{t}}, \tilde{z}_{\tilde{t}})$ satisfy the following FBSDE 
\begin{align}
    \begin{cases}
        d\tilde{x}_{\tilde{t}} &= b(\tilde{t}, \tilde{x}_{\tilde{t}}) d\tilde{t} + \sigma(\tilde{t}, \tilde{x}_{\tilde{t}}) d\tilde{w}_{\tilde{t}}, \\
        d\tilde{y}_{\tilde{t}} &= -g(\tilde{t},\tilde{x}_{\tilde{t}} , \tilde{y}_{\tilde{t}}, \tilde{z}_{\tilde{t}}) d\tilde{t} + \tilde{z}_{\tilde{t}} d\tilde{w}_{\tilde{t}}, \label{transformed FBSDE}\\
        \tilde{y}_{\tilde{T}} &= H(\tilde{x}_{\tilde{T}}).
    \end{cases}
\end{align}
\\
\noindent In this section, we summarize the results of \cite{sym11091153} to compute symmetries of the FBSDE \eqref{FBSDE}. Following the results in \cite{Giuseppe-Gaeta_1999,sym11091153}, we introduce at first order in $\varepsilon$ (or near-identity change of coordinates)
\begin{align}
    \tilde{t}&=t+\varepsilon\tau(t),\label{tilde de t fbsde}\\
     \tilde{x}&=t+\varepsilon\xi(t,x,y),\label{tilde de x fbsde}\\
     \tilde{y}&=t+\varepsilon \eta(t,x,y)\label{tilde de y fbsde}.
\end{align}
\begin{equation}
 d\tilde{\omega}=\left(1+\varepsilon\frac{\tau^{'}(t)}{2}\right)d\omega.\label{tilde de w fbsde oks}
\end{equation}
\begin{thm}\cite[p.~11]{sym11091153}
     We have at first order in $\varepsilon$
\begin{equation}
    \bar{z}=z+\varepsilon \left(\eta_x\sigma+\eta_y z-\frac{\tau_t}{2}z\right).
\end{equation}
\end{thm}
\begin{proof}
    The proof is derived by applying Itô formula to the expansion of the equation \eqref{tilde de y fbsde} and utilizing the fact that the Brownian motion is transformed as described in equation \eqref{tilde de w fbsde oks}.
\end{proof}
\begin{thm}\cite[p.~12]{sym11091153}\label{terminal condition holds}
For terminal condition $Y_T = H(X_T)$, we have that $$\Tilde{Y}_{\Tilde{T}} = H(\Tilde{X}_{\tilde{T}}) \,\,\text{holds if and olny if}\,\, \eta(t,x,H(x))=H_x \xi(t,x,H(x))$$
\end{thm}
\noindent We now have all the necessary tools to state the following theorem.

\begin{thm}{\cite[p.~13]{sym11091153}}\label{DS FBSDE}
The projectable vector field \eqref{Projectable VF} is a symmetry generator of the FBSDE \eqref{FBSDE} if and only if $\tau, \xi, \eta$ satisfy the following determining equations:
    \begin{align}
&(g\tau)_{t}+g_{x}\xi+g_{y}\eta+g_{z}\left(\sigma\eta_{x}+\left(\eta_{y}-\frac{1}{2}\tau_{t}\right)z\right)+\eta_{t}+b\eta_{x}-g\eta_{y}+\frac{1}{2}\left(\sigma^{2}\eta_{xx}
+\eta_{yy}z^{2}+\eta_{xy}\sigma z\right)=0 \label{(gtau)_t FBSDE}\\
&-(b\tau)_{t}-b_{x}\xi+\xi_{t}+b\xi_{x}-g\xi_{y}+\frac{1}{2}\left(\sigma^{2}\xi_{xx}+\xi_{yy}z^{2}+\sigma\xi_{xy}z\right)=0\label{(btau)_t FBSDE}\\
&\frac{1}{2}\sigma\tau_{t}+\sigma_{t}\tau+\sigma_{x}\xi-\sigma\xi_{x}-\xi_{y}z=0 \label{xi_y=0 FBSDE}\\
&\eta(t, x, H(x)) \equiv \partial_x H(x) \xi(t, x, H(x))\label{H}
 \end{align}
\begin{remark}
\begin{itemize}
    \item 
    It is interesting to note that the determining equations in Theorem \ref{DS FBSDE} are deterministic even though they
describe symmetries of a stochastic equations.
\item 
In the general case, when functions  $b(t, x),\,\,\sigma(t, x)$ and $g(t,x,y,z)$ are arbitrary, the determining equations of the FBSDE/PDE have no non-trivial solutions, i.e. there are no symmetries.
\item 
 If $\sigma \ne 0$, $\eta=0$ and $\xi$ depends only on $t,x$, we find the result of  Gaeta \cite{Giuseppe-Gaeta_1999} and Unal \cite{unal2003symmetries}, the determining system of the stochastic differential equations (SDE).
\end{itemize}
  
\end{remark}

\end{thm}
\noindent Since, $\sigma$, $\tau$ and $\xi$ does not depend on $z$, and from \eqref{xi_y=0 FBSDE} we deduce $\xi_{y}=0$ and \\
$\frac{1}{2}\sigma\tau_{t}+\sigma_{t}\tau+\sigma_{x}\xi-\sigma\xi_{x}=0$. Therefore, \eqref{(btau)_t FBSDE} becomes $-(b\tau)_{t}-b_{x}\xi+\xi_{t}+b\xi_{x}+\frac{1}{2}\sigma^{2}\xi_{xx}=0.$ Also condition \eqref{H} becomes $\eta(t, x, H(x)) \equiv \partial_x H(x) \xi(t, x)$. Finally, Theorem \ref{DS FBSDE} becomes:
\begin{thm}\label{X0-FBSDE}
The projectable vector field \eqref{Projectable VF} becomes $v_{_{FBSDE}}=\tau(t) \partial_t+\xi(t, x) \partial_x+\eta(t, x, y) \partial_y$. The final version of the determining system of the FBSDE \eqref{FBSDE}:
\begin{align}
    &(g\tau)_{t}+g_{x}\xi+g_{y}\eta+g_{z}\left(\sigma\eta_{x}+\left(\eta_{y}-\frac{1}{2}\tau_{t}\right)z\right)+\eta_{t}+b\eta_{x}-g\eta_{y}+\frac{1}{2}\left(\sigma^{2}\eta_{xx}
+\eta_{yy}z^{2}+\eta_{xy}\sigma z\right)=0\\
&\xi_{t}+\frac{1}{2}\sigma^{2}\xi_{xx}+b\xi_{x}-b_{x}\xi-(b\tau)_{t}=0\\
&-\sigma\xi_{x}+\sigma_{x}\xi+\frac{1}{2}\sigma\tau_{t}+\sigma_{t}\tau=0\\
&\eta(t, x, H(x)) \equiv \partial_x H(x) \xi(t, x).
\end{align}

\end{thm}

\section{Symmetries of the semi-linear PDE}
We consider the semi-linear PDE with final condition, that is related to the FBSDE \eqref{FBSDE} via proposition \ref{PDE-FBSDE}. We give in the following theorem the equations that determines the structure of the symmetries of the PDE \eqref{SemiLinear-PDE}. Nowadays, the problem of computing Lie symmetries of deterministic PDEs or ODEs is a classic, for example see \cite{OLV86,ovsiannikov2014group,blu,ibragimov1995crc}. Not like symmetries of SDE which is recent and limited.
\begin{defin}
    $$\frac{\partial }{\partial (\sigma(t,x) u_x)}:=\frac{\partial }{\partial z}$$
\end{defin}
\noindent we justify the definition above by proposition \ref{PDE-FBSDE}.
\begin{thm}\label{DS(PDE)}
A vector field
    \begin{equation}
        v_{_{PDE}}=\theta(t,x,u) \partial_t+\gamma(t, x, u) \partial_x+\phi(t, x, u) \partial_u  \label{X_{1}}
    \end{equation}
    is a symmetry generator of the the PDE \eqref{SemiLinear-PDE} if and only if 
    the coefficients $\theta,\,\gamma$ and $\phi$ satisfy the following equations:
    \begin{align*}
    &\Bigg[(g\theta)_{t}+g_{x}\gamma+g_{u}\phi+ \frac{1}{2}\sigma^{2}\phi_{xx}+\left(b+\sigma g_z\right)\phi_{x}+\phi_{t}-g\phi_{u}+\frac{1}{2}\sigma^{2}g\theta_{xx}+\left(b+\sigma g_z\right)g\theta_{x}\\
    &-g^{2}\theta_{u}\Bigg](t,x,u,\sigma u_x)=0\\
    &\Bigg[\frac{1}{2}\sigma^{2}\gamma_{xx}+\left(b+\sigma g_z\right)\gamma_{x}+\gamma_{t}-g\gamma_{u}-\left(b_x+\sigma_x g_z\right)\gamma-\frac{1}{2}\sigma^{2}b\theta_{xx}-\sigma^{2}g\theta_{xu}+\left(b-\sigma g_z\right)g\theta_{u}\\
    &-\left(b+\sigma g_z\right)b\theta_{x}-(b\theta)_{t}-\sigma^{2}\phi_{xu}-\sigma_t g_z\theta-\sigma g_z\phi_{u}\Bigg](t,x,u,\sigma u_x)=0\\
     &\Bigg[\frac{1}{2}\sigma^{4}\theta_{xx}+\frac{1}{2}\sigma^{2}\left(b+\sigma g_z\right)\theta_{x}-\frac{1}{2}\sigma^{2}g\theta_{u}+\frac{1}{2}\sigma^{2}\theta_{t}+\sigma\sigma_{t}\theta-\sigma^{2}\gamma_{x}+\sigma\sigma_{x}\gamma \Bigg] (t,x,u,\sigma u_x)=0\\
     &\Bigg[\frac{1}{2}\sigma^{2}\phi_{uu}-\sigma^{2}\gamma_{xu}-\sigma g_z\gamma_{u}+ \frac{1}{2}\sigma^{2}g\theta_{uu}+b\sigma^{2}\theta_{xu}+b\sigma g_z \theta_{u}\Bigg](t,x,u,\sigma u_x)=0\\
     &\Bigg[\frac{1}{2}\sigma^{3} g_z\theta_{u}+\frac{1}{2}\sigma^{4}\theta_{xu}-\sigma^{2}\gamma_{u}\Bigg](t,x,u,\sigma u_x)=0\\
     &\frac{1}{2}\sigma^{2}\left(-\gamma_{uu}+b\theta_{uu}\right)=0,\quad\sigma^{2}\theta_{x}=0,\quad\sigma^{2}\theta_{u}=0\\
     &\phi(t, x, H(x)) \equiv \partial_x H(x) \gamma(t, x, H(x))\,\,(\text{so the terminal condition holds}).
    \end{align*}
\end{thm}
\begin{proof}
We give only essential steps for the proof. Recall the PDE \eqref{SemiLinear-PDE},
\begin{align*}
    \begin{cases}
        &\partial_t u(t, x)+b(t,x)\partial_x u(t, x)+\frac{1}{2}\sigma^2(t,x)\partial^{2}_{xx} u(t, x)+g\left(t, x, u(t, x), \partial_x u(t, x) \sigma(t, x)\right)=0, \\
        &u(T, x)=H(x), \quad x \in \mathbb{R}.
    \end{cases}
\end{align*}
    We define the following spaces, which will be used throughout the proof:\\
    $X$: the set of independent variables of the PDE \eqref{PDE-FBSDE}, with elements denoted by $x$; typically, we take $X\subset\mathbb{R}^p$, where $p$ is number of independent variables.\\
    $U$: the set of dependent variables of the PDE \eqref{PDE-FBSDE}, with elements denoted by $u$, typically, we take $U\subset\mathbb{R}^q$, where $q$ is number of dependent variables.\\
    $U^{(n)}$: is the set of all derivatives of $u$ with respect to independent variables, up to order $n$, we can express $U^{(n)}$ as the cartesian product $U^{(n)}=U\times U_{1}\times...\times U_{n}$, where $U_k$ denotes the space of $k$-th order derivatives. This space is referred to as the $n$-th jet space of $U$.
In this case, we have 2 independents variables: $x$ and $t$, along with one dependent variable $u$. The semi linear PDE is of second order, and we can identify it as a submanifold in the jet space $X\times U^{(2)}$ determined by the vanishing of
$$\Delta(t,x,u^{(2)})=\partial_t u(t, x)+b(t,x)\partial_x u(t, x)+\frac{1}{2}\sigma^2(t,x)\partial^{2}_{xx} u(t, x)+g\left(t, x, u(t, x), \partial_x u(t, x) \sigma(t, x)\right).$$
Let 
$$v_{_{PDE}}=\theta(t,x,u) \partial_t+\gamma(t, x, u) \partial_x+\phi(t, x, u) \partial_u $$
 be a vector field on $X\times U$. We want to determine all conditions on coefficients $\theta$, $\gamma$ and $\phi$ so that the corresponding one-parameter group $exp(\varepsilon v)$ is a symmetry group of \eqref{SemiLinear-PDE}.
According to Lie's Theorem in \cite[p.~104]{OLV86}, we need to know the second order prolonged vector field $pr^{2}v$ defined in the prolonged jet space $X\times U^{(2)}$. Using the general prolongation formula in  \cite[p.~110]{OLV86} , we can write the prolonged vector field $v_{_{PDE}}$ as 
$$pr^{2}v_{_{PDE}}=v_{_{PDE}}+\phi^{x}\frac{\partial}{\partial u_{x}}+\phi^{t}\frac{\partial}{\partial u_{t}}+\phi^{xx}\frac{\partial}{\partial u_{xx}}+\phi^{xt}\frac{\partial}{\partial u_{xt}}+\phi^{tt}\frac{\partial}{\partial u_{tt}},
$$ 
whose coefficients were calculated in example 2.38 in \cite[p.~114]{OLV86}.\\
Applying Lie's Theorem in \cite[p.~104]{OLV86} we have that $v_{_{PDE}}$ is a symmetry generator of \eqref{SemiLinear-PDE} if and only if $$pr^{2}v_{_{PDE}}\left(\Delta(t,x,u^{(2)})\right)=0 \,\, \text{such that }\,\, \Delta(t,x,u^{(2)})=0$$ 
We develop $pr^{2}v_{_{PDE}}\left(\Delta(t,x,u^{(2)})\right)$ and replace $u_t$ by \\$-\left(b(t,x)\partial_x u(t, x)+\frac{1}{2}\sigma^2(t,x)\partial^{2}_{xx} u(t, x)+g\left(t, x, u(t, x), \partial_x u(t, x) \sigma(t, x)\right)\right)$ whenever it occurs. By equating the coefficients of the various monomials in the first and second order partial derivatives of $u$, finally we derive the determining equations (stated in the theorem above) for the symmetry groups of the semi linear PDE \eqref{SemiLinear-PDE}.
\end{proof}

\begin{remark}
\begin{itemize}
    \item  
    If $g=0$, we find the determining system of the Kolmogorov backward equation, found by Kozlov \cite{kozlov2021symmetries}.
    \item 
  If $g$ is linear in $u$ which corresponds to the Feynman-Kac formula, numerous authors have examined the symmetries of the corresponding PDE:
   \begin{align*}
      u_t&=\sigma x^{\gamma}u_{xx}+f(x)u_x-\mu x^{r}u,\\
       u_t&=\sigma x^{\gamma}u_{xx}+f(x)u_x+g(x)u.
  \end{align*}
  For instance, see \cite{craddock2020lie,craddock2010equivalence,craddock2009calculation,craddock2007lie}.
  \item 
   If $\sigma(t,x)=constante\ne 0$, we find the result of \cite[p.~2]{lescot2019symmetries}.
\end{itemize} 
\end{remark}
\begin{corollary}
    If $\sigma\neq0$, every symmetry of the PDE \eqref{SemiLinear-PDE} is projectable.
\end{corollary}

\noindent We are going to simplify the determining system in Theorem $\eqref{DS(PDE)}$ as follows.
\begin{corollary}
If $\sigma(t,x) \ne 0$, the determining system of the PDE \eqref{SemiLinear-PDE} becomes: 
     \begin{align*}
         &\Bigg[\phi_{t}+\frac{1}{2}\sigma^{2}\phi_{xx}+\left(b+\sigma g_z\right)\phi_{x}+(g\theta)_{t}+g_{x}\gamma+g_{u}\phi-g\phi_{u}\Bigg](t,x,u,\sigma u_x)=0,\\      
          &\Bigg[-\gamma_{t}-\frac{1}{2}\sigma^{2}\gamma_{xx}-\left(b+\sigma g_z\right)\gamma_{x}+(b_{x}+\sigma_x g_z)\gamma+\sigma_t g_z\theta+(b\theta)_{t}+
\sigma g_z\phi_u\\
&+\sigma^{2}\phi_{xu}\Bigg](t,x,u,\sigma u_x)=0,\\      
        &-\sigma\gamma_{x}+\sigma_{x}\gamma+\frac{1}{2}\sigma\theta_{t}+\sigma_{t}\theta=0,\\
        &\phi_{uu}=0,\,\,\theta_{x}=\theta_{u}=\gamma_{u}=0,\\
        &\phi(t, x, H(x)) \equiv \partial_x H(x) \gamma(t, x)\,\,(\text{so the terminal condition holds}).
    \end{align*}
\end{corollary}
\noindent  Using proposition \ref{PDE-FBSDE}, the determining system of the PDE \eqref{SemiLinear-PDE} becomes 
\begin{align}
 &\phi_{t}+\frac{1}{2}\sigma^{2}\phi_{xx}+(b+\sigma g_z)\phi_{x}+(g\theta)_{t}+g_{x}\gamma+g_{y}\phi-g\phi_{y}=0 \label{phi_{t}}\\      
          &-\gamma_{t}-\frac{1}{2}\sigma^{2}\gamma_{xx}-(b+\sigma g_z)\gamma_{x}+(b_{x}+\sigma_x g_z)\gamma+\sigma_t g_z\theta+(b\theta)_{t}+
\sigma g_{z}\phi_y+\sigma^{2}\phi_{xy}=0 \label{gamma_{t}}\\      
       & -\sigma\gamma_{x}+\sigma_{x}\gamma+\frac{1}{2}\sigma\theta_{t}+\sigma_{t}\theta=0 \label{sigmagamma_{x}}\\
        &\phi_{yy}=0 \label{phi_yy=0}\\
        &\phi(t, x, H(x)) \equiv \partial_x H(x) \gamma(t, x) \label{final condition}.
\end{align}
\noindent Condition \eqref{phi_yy=0} gives that 
$\phi(t,x,y)=\delta(t,x) y+\beta(t,x)$ for certain functions $\delta$ and $\beta$.
So, the final structure of the determining system of the PDE \eqref{SemiLinear-PDE} become

\begin{align}
&\phi(t,x,y)=\delta(t,x) y+\beta(t,x) \\
 &\Bigg[\frac{1}{2}\sigma^{2}(\beta_{xx}+\delta_{xx}y)+(b+\sigma g_z)(\beta_x+\delta_x y)+(\beta_{t}+\delta_{t}y)+g_{y}(\beta+\delta y)-g\delta\\
 \nonumber &+g_{x}\gamma+(g\theta)_{t}\Bigg](t,x,u,\sigma u_x)=0\\     
 &\Bigg[-\frac{1}{2}\sigma^{2}\gamma_{xx}-\gamma_{t}-(b+\sigma g_z)\gamma_{x}+(b_{x}-\sigma_x g_z)\gamma+(b\theta)_{t}+\sigma^{2}\delta_{x}+\sigma g_z \delta\\
 \nonumber&+\sigma_t g_z \theta\Bigg](t,x,u,\sigma u_x)=0 \\      
        &\frac{1}{2}\sigma\theta_{t}+\sigma_{t}\theta-\sigma\gamma_{x}+\sigma_{x}\gamma=0 \\
        &\beta(t,x)=\partial_{x}H(x)\gamma(t,x)-\delta(t,x)H(x).
\end{align}
\noindent Therefore, for $\sigma\neq 0$, the symmetry generator of the PDE \eqref{SemiLinear-PDE} is of the following form 
\begin{equation}
    v_{_{PDE}}=\theta(t) \partial_t+\gamma(t,x) \partial_x+\left[\delta(t,x)y+\partial_x H(x)\gamma(t,x)-\delta(t,x)H(x)\right]\partial_y
\end{equation}
for every given $H$.\\
\section{Symmetry of FBSDE versus symmetry of semi linear PDE}
Now, we are going to analyse and determine the relation between symmetries of FBSDE and  symmetries of the associated non-linear PDE. Inspired by the works of Gaeta \cite{Giuseppe-Gaeta_1999} and Kozlov \cite{kozlov2021symmetries}, one might expect that the symmetries of FBSDEs are contained within the symmetries of the nonlinear PDE, or vice versa, or that they share some common symmetries. However, we previously observed that for arbitrary $b,\,\sigma$, and $g$, the determining system of the FBSDE \eqref{FBSDE} and the PDE has no solutions. Therefore, there is no basis for comparing the symmetries of FBSDEs and PDEs. However, for certain specific (less general) forms of $b$ or $\sigma$ or $g$, we hope to gain more informations about $\tau,\,\theta,\,\xi,\,\gamma,\,\eta$, and $\phi$. This would enable us to compare the symmetries of FBSDEs \eqref{FBSDE} and semi linear PDEs \eqref{SemiLinear-PDE}.
\begin{defin}
    We call
    \begin{align*}
        &\mathcal{L}_{ FBSDE}: \,\,\text{Symmetry space of the FBSDE}\,\, \eqref{FBSDE},\\ &\mathcal{L}_{PDE}: \,\,\text{Symmetry algebra of the semi linear PDE}\,\, \eqref{SemiLinear-PDE},\\
        &\widehat{\mathcal{L}}_{PDE}:=\left\{ v_{_{PDE}}\in\mathcal{L}_{PDE}/\delta_x=0\right\}.
    \end{align*}
    \end{defin}
\begin{remark}
    \begin{itemize}
    \item 
    $\widehat{\mathcal{L}}_{PDE}$ is a subalgebra of $\mathcal{L}_{PDE}$.
        \item 
        $\mathcal{L}_{ FBSDE}$ is not in general a Lie algebra. For a counter example, see \cite[p.~12]{ZHANG2021105527}.
    \end{itemize}
\end{remark}   
\noindent For $g$ independent of the variable $z$, and since $\tau,\,\theta,\,\xi,\,\gamma,\,\eta$, and $\phi$ does not depend on $z$ then we can simplify the determining equations for both FBSDE \eqref{FBSDE} and the PDE \eqref{SemiLinear-PDE}. we obtain more informations about $\eta$ and $\phi$. For that we have the following theorem.
\begin{thm}
       If $g$ independent of $z$ then $\mathcal{L}_{ FBSDE}=\widehat{\mathcal{L}}_{PDE}\subseteq \mathcal{L}_{PDE}$.
\end{thm}
\noindent This will indeed appear in the next example.
\begin{proof}
\noindent We will now proceed to prove the stated theorem. Using theorem \ref{X0-FBSDE}, and the fact that $\tau,\, \xi,\,\eta$ and $g$ does not depend on $z$, we find that the determining system of the FBSDE become:

    \begin{align} 
&\frac{1}{2}\sigma^{2}\eta_{xx}+b\eta_{x}+\eta_{t}+ (g\tau)_{t}+g_{x}\xi-g\eta_{y}+g_{y}\eta=0\\
&\eta_{yy}=0\label{eta_yy=0(g<>z)}\\
&\sigma\eta_{xy}=0\label{eta_xy=0(g<>z)}\\
&-(b\tau)_{t}-b_{x}\xi+\xi_{t}+b\xi_{x}+\frac{1}{2}\sigma^{2}\xi_{xx}=0\\
&\frac{1}{2}\sigma\tau_{t}+\sigma_{t}\tau+\sigma_{x}\xi-\sigma\xi_{x}=0\\
&\eta(t, x, H(x)) \equiv \partial_x H(x) \xi(t,x)
 \end{align}

\noindent Equation \eqref{eta_yy=0(g<>z)}, gives that $\eta$ is affine on $y$, $\eta(t,x,y)=\lambda(t,x)y+\mu(t,x)$
for certain function $\lambda$ and $\mu$.  Using equation \eqref{eta_xy=0(g<>z)}, and the fact that $\sigma\ne 0$, we have that $\lambda$ does not depend on $x$. $\eta$ becomes $\eta(t,x,y)=\lambda(t)y+\mu(t,x)$.  Therefore, the final structure of the determining system of the FBSDE is 

    \begin{align} 
    &\eta(x,y,t)=\lambda(t)y+\mu(t,x)\\
&\frac{1}{2}\sigma^{2}\mu_{xx}+b\mu_{x}+\mu_{t}+g_{y}\mu+\lambda_{t}y+\lambda g_{y}y-g\lambda+g_{x}\xi+(g\tau)_{t}=0\\
&\frac{1}{2}\sigma^{2}\xi_{xx}+b\xi_{x}+\xi_{t}-b_{x}\xi-(b\tau)_{t}=0\\
&\frac{1}{2}\sigma\tau_{t}+\sigma_{t}\tau-\sigma\xi_{x}+\sigma_{x}\xi=0\\
&\mu(t,x)=\partial_{x}H(x)\xi(t,x)-\lambda(t)H(x)
 \end{align}
\noindent Finally, the symmetry generator of the FBSDE is of the following form 
\begin{equation}
    v_{_{FBSDE}}=\tau(t) \partial_t+\xi(t,x) \partial_x+\left[\lambda(t)y+\partial_x H(x)\xi(t,x)-\lambda(t)H(x)\right]\partial_y
\end{equation}
for every given $H$.

\noindent It is observed that if $\delta$ is independent of $x$, then,  $\theta,\gamma,\delta,\beta$ and $\tau,\xi,\lambda,\mu$ satisfy the determining systems of the FBSDE and the semi-linear PDE, means they share the same symmetries. Otherwise, the semi-linear PDE exhibits more symmetries (including those of FBSDE). Finally, $\mathcal{L}_{ FBSDE}\subseteq \mathcal{L}_{PDE}$ when $g$ is independent of $z$.
\end{proof}
\section{Application}
In this section we treat some simple cases to illustrate our results and to check the correspondence between symmetries of a FBSDE \eqref{FBSDE} and the associated semi-linear PDE \eqref{SemiLinear-PDE}.
\subsection{Heat equation}
First, we consider the well know heat equation with final condition

  $$
\left\{\begin{array}{l}
u_{t}+\frac{1}{2}u_{xx}=0,\\
u(T,x)=H(x).
\end{array}\right.
$$

\noindent We can associate to this PDE with final condition the following FBSDE 
  $$
\left\{\begin{array}{l}
 dX_t =dW_t, \\
        dY_t =Z_t dW_t, \\
        Y_T = H(X_T).
\end{array}\right.
$$
\noindent Here $g(t,x,y,z)=0$, $b(t,x)=0$ and $\sigma(t,x)=1$.\\
The determining system of the heat equation with final condition is 




\begin{align}
 (\delta_{t}+\frac{1}{2}\delta_{xx})y+(\beta_{t}+\frac{1}{2}\beta_{xx})&=0,\label{delta&beta}\\      
          -(\beta_{t}+\frac{1}{2}\beta_{xx})+\delta_{x}&=0,\label{delta_{x}}\\      
        \beta(t,x)&=\frac{1}{2}\theta_{t}x+h(t),\label{beta}\\
        \phi(t,x,y)&=\delta(t,x)y+\beta(t,x),\\
        \delta(t,x)H(x)+\beta(t,x)&=\partial_x H(x) \gamma(t, x)\label{final condition*}.
\end{align}
\noindent Equation \eqref{delta&beta}, gives that $\delta$ and $\beta$ are solutions of the heat equation.
Replacing $\beta$ in equation \eqref{delta_{x}}, we find that $\delta(t,x)=\frac{1}{4}\theta_{tt}x^{2}+h_{t}x+s(t)$, for certain function $s$. Equation \eqref{final condition*} becomes 
$$
    \beta(t,x)=\partial_x H(x) \left(\frac{1}{2}\theta_{t}x+h(t)\right)-\left(\frac{1}{4}\theta_{tt}x^{2}+h_{t}x+s(t)\right)H(x)
$$
Since $\delta$ and $\beta$ are solutions of the heat equation, then 
\begin{align}
   &\frac{1}{4}\theta_{ttt}x^{2}+h_{tt}x+s_{t}+\frac{1}{4}\theta_{tt}=0\label{theta_{ttt}},\\
   &\left(\frac{1}{2}\theta_{t}x+h(t)\right) H_{xxx}+\left(-\frac{1}{4}\theta_{tt}x^{2}-h_{t}x-s(t)+\theta_{t}\right) H_{xx}=0.\label{H_{xxx}}
   \end{align}
So, the equations that determines symmetries of the heat equation with finale condition are :
\begin{align}
     &\phi(t,x,y)=\delta(t,x)y+\beta(t,x),\\
    &\delta(t,x)=\frac{1}{4}\theta_{tt}x^{2}+h_{t}x+s(t),\\
    &\beta(t,x)=\partial_x H(x) \left(\frac{1}{2}\theta_{t}x+h(t)\right)-\left(\frac{1}{4}\theta_{tt}x^{2}+h_{t}x+s(t)\right)H(x),\\
    &\frac{1}{4}\theta_{ttt}x^{2}+h_{tt}x+s_{t}+\frac{1}{4}\theta_{tt}=0\label{theta_{ttt}},\\
    &\left(\frac{1}{2}\theta_{t}x+h(t)\right) H_{xxx}+\left(-\frac{1}{4}\theta_{tt}x^{2}-h_{t}x-s(t)+\tau_{t}\right) H_{xx}=0.
\end{align}
\noindent Now, Let's find the determining system of the FBSDE associated to the heat equation with final condition. Using the determining equations in theorem (3.4), and the fact that $g$, $\tau$, $\xi$ and $\eta$ does not depend on z, then the determinig equation become 

    \begin{align} 
    &\eta(x,y,t)=\lambda(t)y+\mu(t,x)\\
&\frac{1}{2}\mu_{xx}+\mu_{t}+\lambda_{t}y=0\label{lambda_t*y}\\
&\frac{1}{2}\xi_{xx}+\xi_{t}=0\label{1/2xi_xx+xi_t}\\
&\frac{1}{2}\tau_{t}-\xi_{x}=0\label{tau_t-xi_x}\\
&\mu(t,x)=\partial_{x}H(x)\xi(t,x)-\lambda(t)H(x)\label{mu=H_x*xi-lambda*H}
 \end{align}
\noindent Using equation \eqref{lambda_t*y}, and since $\mu$ and $\lambda$ does not depend on $y$ then 
\begin{align}
    \lambda_{t}&=0,\\
    \frac{1}{2}\mu_{xx}+\mu_{t}&=0\label{1/2*mu_xx+mu_t}.
\end{align}
Equation \eqref{tau_t-xi_x} gives $\xi(t,x)=\frac{1}{2}\tau_{t}x+l(t)$ for certain function $l$. Replacing $\xi$ in equation \eqref{1/2xi_xx+xi_t} gives $\frac{1}{2}\tau_{tt}x+l_{t}=0$, since $\tau$ and $l$ depends only on $t$, then $\tau_{tt}=0$  and $h_{t}=0$. Replacing $\mu$ of  equation \eqref{mu=H_x*xi-lambda*H} in equation \eqref{1/2*mu_xx+mu_t}, we find that 

$$\frac{1}{2}\left[(\frac{1}{2}\tau_{t}x+h)H_{xxx}+(\tau_{t}-c_{1})H_{xx}\right]+\left[\frac{1}{2}\tau_{tt}x+h_{t}\right]H_{x}=0.$$

\noindent It provides the final structure of the determining equations of the FBSDE,
    \begin{align} 
    \eta(x,y,t)=\lambda(t) y+\mu(t,x)\\
    \xi(t,x)=\frac{1}{2}\tau_{t}x+h(t)\\
    \lambda_{t}=0\\
    \tau_{tt}=0\\
    h_{t}=0\\
    (\frac{1}{2}\tau_{t}x+h)H_{xxx}+(\tau_{t}-c_{1})H_{xx}=0\\
\mu(t,x)=\partial_{x}H(x)(\frac{1}{2}\tau_{t}x+h(t))-c_{1}H(x)
 \end{align}
\noindent The FBSDE and the PDE are defined for any given $H$. Let us present some examples of $H$ to clearly illustrate how the symmetries of the FBSDE are included in the symmetries of the PDE. 
\begin{enumerate}
    \item \textbf{If H}$(x)=x$\\
    Solving the determining equation of the FBSDE yields:
    \begin{align} 
    \eta(x,y,t)=c_{1}y+(\frac{1}{2}c_{2}-c_{1})x+c_{4}\\
    \xi(t,x)=\frac{1}{2}c_{2}x+c_{4}\\
    \tau(t)=c_{2}t+c_{3}
 \end{align}
Thus, the symmetries of the FBSDE are represented by the following vector fields:
\begin{align} 
    v_{1}=(y-x)\partial_{y}\\
    v_{2}=t\partial_{t}+\frac{1}{2}x\partial_{x}+\frac{1}{2}x\partial_{y}\\
    v_{3}=\partial_{t}\\
    v_{4}=\partial_{x}+\partial_{y}
 \end{align}
Computing the Lie bracket:
\begin{longtable}{|c|c|c|c|c|} 
\hline
$[v_{i},v_{j}]$ & $v_{1}$ & $v_{2}$&$v_{3}$&$v_{4}$ \endhead
\hline
$v_{1}$&0&$0$&$0$&0\\ 
$v_{2}$&$0$&0&$0$&0\\
$v_{3}$&$0$&$0$&0&0\\
$v_{4}$&0&0&0&0\\
\hline
\end{longtable}
From the commutator table, we see that the symmetries of the FBSDE forms an Abelian Lie algebra of dimension 4. Solving the determining equation of the PDE gives:
\begin{align}
     &\phi(t,x,y)=\left(\left(\frac{1}{4}x^{2}-\frac{1}{4}t\right)c_{1}-\frac{1}{4}c_{2}+c_{4}x+c_{6}\right)(y-x)+\frac{1}{2}txc_{1}+\frac{1}{2}xc_{2}+c_{4}t+c_{5},\\
   &\gamma(t,x)=\frac{1}{2}txc_{1}+\frac{1}{2}xc_{2}+c_{4}t+c_{5},\\
   &\theta(t)=\frac{1}{2}t^{2}c_{1}+c_{2}t+c_{3}.
\end{align}
it gives the symmetries of the heat equation:

\begin{align}
&w_{1}=\frac{1}{2}t^{2}\partial_{t}+\frac{1}{2}tx\partial_{x}+\left(\left(\frac{1}{4}x^{2}-\frac{1}{4}t\right)(y-x)+\frac{1}{2}tx\right)\partial_{y},\\
&w_{2}=t\partial_{t}+\frac{1}{2}x\partial_{x}+\left(\frac{1}{2}x-\frac{1}{4}(y-x)\right)\partial_{y},\\
 &w_{3}=\partial_{t},\\
        &w_{4}=t\partial_{x}+(x(y-x)+t)\partial_{y},\\
         &w_{5}=\partial_{x}+\partial_{y},\\
          &w_{6}=(y-x)\partial_{y}.
\end{align}
computing the Lie bracket:
\begin{longtable}{|c|c|c|c|c|c|c|} 
\hline
$[w_{i},w_{j}]$ & $w_{1}$ & $w_{2}$&$w_{3}$&$w_{4}$&$w_{5}$&$w_{6}$ \endhead
\hline
$w_{1}$&0&$-w_{1}$&$-w_{2}$&0&$-\frac{1}{2}w_{4}$&0\\ 
$w_{2}$&$w_{1}$&0&$-w_{3}$&$\frac{1}{2}w_{4}$&$-\frac{1}{2}w_{5}$&0\\
$w_{3}$&$w_{2}$&$w_{3}$&0&$w_{5}$&$0$&0\\
$w_{4}$&0&$-\frac{1}{2}w_{4}$&$-w_{5}$&0&$-w_{6}$&0\\
$w_{5}$&$\frac{1}{2}w_{4}$&$\frac{1}{2}w_{5}$&$0$&$w_{6}$&0&0\\
$w_{6}$&0&0&0&0&0&0\\
\hline
\end{longtable}
\begin{remark}
\begin{itemize}
    \item The Lie algebra of the Heat equation with terminal condition called $\mathcal{G}$ is isomorphic to $\mathfrak{sl}_2\ltimes \mathfrak{H}_3\cong \mathcal{H}_{0,0}$, where $\mathcal{H}_{0,0}$ is the algebra of heat equation found by \cite{lescot2014solving}. 
    \item 
   By imposing the terminal condition $u(T,x)=x$ (with $H(x)=x$), the infinite dimensional Lie Algebra spanned by $v_{\mu}=\mu(t,x)\frac{\partial}{\partial u}$, where $\mu$ is an arbitrary solution of the heat equation $u_t=u_{xx}$ described in Olver \cite{OLV86}, is lost. Consequently, as previous stated, we obtain only a finite-dimensional Lie algebra.
\end{itemize}   
\end{remark}
Finally, we see that, the symmetries of the FBSDE are included in the symmetries of the PDE in the case of $H(x)=x$.

 \item \textbf{If H}$(x)=x^{2}$\\
The solution of the determining equations of the PDE is:

\begin{align}
&\theta(t)=c_{1}t+c_{2},\\
&\gamma(t,x)=\frac{1}{2}xc_{1}+c_{3},\\
&\phi(t,x,y)=c_{1}y+2xc_{3}.
\end{align}
The symmetry algebra is spanned by 3 symmetries:
\begin{align}
&w_{1}=t\partial_{t}+\frac{1}{2}x\partial_{x}+y\partial_{y},\\
&w_{2}=\partial_{t},\\
&w_{3}=\partial_{x}+2x\partial_{y}.
\end{align}
Solving now, the determining equations of the FBSDE gives the following symmetries:
\begin{align}
&v_{1}=t\partial_{t}+\frac{1}{2}x\partial_{x}+y\partial_{y},\\
&v_{2}=\partial_{t},\\
&v_{3}=\partial_{x}+2x\partial_{y}.
\end{align}
We remark here that both of the determing system of the PDE and FBSDE have the same symmetries.\\
Computing the Lie bracket:
\begin{longtable}{|c|c|c|c|c|} 
\hline
$[v_{i},v_{j}]$ & $v_{1}$ & $v_{2}$&$v_{3}$ \endhead
\hline
$v_{1}$&0&$-v_{2}$&$-\frac{1}{2}v_{3}$\\ 
$v_{2}$&$v_{2}$&0&$0$\\
$v_{3}$&$\frac{1}{2}v_{3}$&$0$&0\\
\hline
\end{longtable}
This Lie algebra is isomorphic to $L\left(3,2,\frac{1}{2}\right)$ using the notations of Bowers \cite{Bowers}.
  \item \textbf{If H}$(x)=x^{n},\,\,n\geq 3$\\
\end{enumerate}
Solving the determining equations of the PDE, gives that 
    \begin{align}
&\theta(t)=c_{1}t+c_{2},\\
&\gamma(t,x)=\frac{1}{2}xc_{1}+c_{3},\\
&\phi(t,x,y)=\frac{n}{2}yc_{1}+2x^{n-1}c_{3}.
\end{align}
 \begin{align}
\text{ Symmetries are:}\quad&w_{1}=t\partial_{t}+\frac{1}{2}x\partial_{x}+\frac{n}{2}y\partial_{y},\\
&w_{2}=\partial_{t},\\
&w_{3}=\partial_{x}+nx^{n-1}\partial_{y}.
\end{align}
\begin{longtable}{|c|c|c|c|c|} 
\hline
$[w_{i},w_{j}]$ & $w_{1}$ & $w_{2}$&$w_{3}$ \endhead
\hline
$w_{1}$&0&$-w_{2}$&$-\frac{1}{2}w_{3}$\\ 
$w_{2}$&$w_{2}$&0&$0$\\
$w_{3}$&$\frac{1}{2}w_{3}$&$0$&0\\
\hline
\end{longtable}
\noindent After solving the determinign system of FBSDE we find that symmetries are:
\begin{align}
&v_{1}=t\partial_{t}+\frac{1}{2}x\partial_{x}+\frac{n}{2}y\partial_{y},\\
&v_{2}=\partial_{t}.
\end{align}
\begin{longtable}{|c|c|c|c|c|} 
\hline
$[v_{i},v_{j}]$ & $v_{1}$ & $v_{2}$ \endhead
\hline
$v_{1}$&0&$-v_{2}$\\ 
$v_{2}$&$v_{2}$&0\\
\hline
\end{longtable}
\noindent This Lie algebra is isomorphic to $L\left(2,1\right)$ using the notations of Bowers \cite{Bowers}.\\
Finally, for $n\geq 3$, symmetries of FBSDE are included in symmetries of the PDE.\\
However, if $g$ depends on $z$, it becomes generally unclear how to establish a direct comparison between the symmetries of the FBSDE \eqref{FBSDE} and those of the semi-linear PDE \eqref{PDE-FBSDE}. Nonetheless, we noticed that in certain specific cases where $g$ depends on $z$, it is possible to transform the FBSDE \eqref{FBSDE} into a new FBSDE in which the modified $g$ no longer depends on $z$.
\subsubsection{Application of Girsanov transformation}
Consider the FBSDE \eqref{FBSDE} with $g(t,x,y,z)=q(t,x,y)+\alpha z$
\begin{align}
    \begin{cases}
        dX_t &= b(t, X_t) dt + \sigma(t, X_t) dW_t, \label{FBSDE-girsanov}\\
        dY_t &= -(q(t, X_t, Y_t)+\alpha Z_t) dt + Z_t dW_t, \\
        Y_T &= H(X_T).
    \end{cases}
\end{align}
where $\alpha$ is a constant. Focusing on the BSDE $dY_t= -(g(t, X_t, Y_t)+\alpha Z_t) dt + Z_t dW_t$.\\
Rewrite it in the integrale form $$Y_t=Y_T +\int_{t}^{T}(q(s, X_s, Y_s)+\alpha Z_s) ds - \int_{t}^{T}Z_s dW_s$$
Equivalent to say $$Y_t=Y_T +\int_{t}^{T}q(s, X_s, Y_s) ds - \int_{t}^{T}Z_s d(W_s-\alpha s)$$
Set $\Bar{W_s}=W_s-\alpha s$\\
Using Girsanov theorem \cite{oksendal2013stochastic} $\Bar{W_s}$ is a also Brownian motion under a new probability $\mathbb{Q}$, and we have $\frac{d\mathbb{Q}}{d\mathbb{P}}_{|_{\mathcal{F}_{t}}}=\exp{\left(\alpha W_{t}-\frac{\alpha^{2}}{2}t\right)}$ due to Girsanov theorem.\\
The BSDE in \eqref{FBSDE-girsanov} becomes $$Y_t=H(X_T) +\int_{t}^{T}q(s, X_s, Y_s) ds - \int_{t}^{T}Z_s d\Bar{W_{s}}.$$
under the new probability $\mathbb{Q}$.\\
Consider now the forward SDE in \eqref{FBSDE-girsanov} $$X_t=X_0 +\int_{0}^{t}b(s, X_s) ds + \int_{0}^{t}\sigma(s,X_s) dW_{s},$$
under the probability $\mathbb{P}$. The SDE becomes 
$$X_t=X_0 +\int_{0}^{t}(b(s, X_s)+\alpha \sigma(s,X_s)) ds + \int_{0}^{t}\sigma(s,X_s) d(W_{s}-\alpha s).$$
Thus, $$X_t=X_0 +\int_{0}^{t}\Bar{b}(s, X_s) ds + \int_{0}^{t}\sigma(s,X_s) d\Bar{W_{s}},$$
under the probability $\mathbb{Q}$, where $\bar{b}(t,x)=b(t,x)+\alpha \sigma(t,x).$
Finally, the FBSDE \eqref{FBSDE-girsanov} becomes 
\begin{align*}
    \begin{cases}
        dX_t &= \Bar{b}(s, X_s) dt + \sigma(t, X_t) d\Bar{W_{s}}, \\
        dY_t &= -q(t, X_t, Y_t) dt + Z_t d\Bar{W_{s}}, \\
        Y_T &= H(X_T).
    \end{cases}
\end{align*}
under the probability $\mathbb{Q}$. Therefore we return to our previous case $g$ independant of $z$.

Inspired by Barrieu and El Karoui in \cite{10.1214/12-AOP743}, who examined the existence of solutions for BSDEs with a generator $g$
that is quadratic in $z$, we are going to analyse the FBSDE \eqref{FBSDE}, where $g$ is quadratic on $z$, more precisely $g$ is of the form $g(t,x,y,z)=q(t,x,y)+r(y)z^{2}$.
\subsubsection{Quadratic case}
Consider the following FBSDE:
\begin{align}
    \begin{cases}
        dX_t &= b(s, X_s) dt + \sigma(t, X_t) dW_{s}, \\
        dY_t &= -(q(t, X_t, Y_t)+r(Y_t)Z_t^{2}) dt + Z_t dW_{s},\label{FBSDE in case of g=q+rZ^2} \\
        Y_T &= H(X_T).
    \end{cases}
\end{align}
We are going to transform the FBSDE \eqref{FBSDE in case of g=q+rZ^2} to the case where the generator $g$ depends only on $t,x,y$. Consider a transformation of the process $Y_t$, called $\bar{Y_t}=\rho(Y_t)$, that we will determine later. Using Ito formula
\begin{align}
    \rho(Y_T)&=\rho(Y_t)+\int_t^T\rho^{'}(Y_s)dY_s+\frac{1}{2}\int_t^T\rho^{''}(Y_s)Z_s^2ds\\
    &=\rho(Y_t)+\int_t^T\rho^{'}(Y_s)[-(q(s, X_s, Y_s)+r(Y_s)Z_s^2) ds + Z_sdW_s]+\frac{1}{2}\int_t^T\rho^{''}(Y_s)Z_s^2ds\\
    &=\rho(Y_t)+\int_t^T-\left[(\rho^{'}q)(s,X_s,Y_s)+(\rho^{'}r)(Y_s)Z_s^2\right] ds +\int_t^T Z_sdW_s+\frac{1}{2}\int_t^T\rho^{''}(Y_s)Z_s^2ds\\
    &=\rho(Y_t)+\int_t^T-(\rho^{'}q)(s,X_s,Y_s) ds +\int_t^T Z_sdW_s+\int_t^T\left[\frac{1}{2}\rho^{''}-(\rho^{'}r)\right] (Y_s)Z_s^2ds
\end{align}
We choose $\rho$ such that 
\begin{equation}
    \frac{1}{2}\rho^{''}-\rho^{'}r=0.\label{ODE(rho)}
\end{equation}
\noindent Therefore, 
\begin{equation}
    \rho(Y_T)=\rho(Y_t)+\int_t^T-(\rho^{'}q)(s,X_s,Y_s) ds +\int_t^T Z_sdW_s
\end{equation}
so, 
\begin{equation}
    \bar{Y_t}=\bar{Y_T}+\int_t^T\rho^{'}\circ\rho^{-1}(\bar{Y_s})q(s,X_s,\rho^{-1}(\bar{Y_s})) ds -\int_t^T Z_sdW_s
\end{equation}
Finally, 
\begin{equation}
    \bar{Y_t}=(\rho\circ H)(X_T)+\int_t^T\bar{q}(s,X_s,\bar{Y_s}) ds -\int_t^T Z_sdW_s
\end{equation}
where, $\bar{q}(s,x,\bar{y})=\rho^{'}\circ\rho^{-1}(\bar{y})q(s,x,\rho^{-1}(\bar{y}))$\\
$$\text{We have then, }\qquad(X,Y,Z;Y_T=H(X_T))\qquad \xrightarrow{\rho} \qquad(X,\rho(Y),Z,\bar{Y_T}=\rho\circ H(X_T))$$
It remains to find $\rho$ by solving the ODE \eqref{ODE(rho)}, solving we find 
\begin{equation}
\rho(p)=\rho(0)+\int_0^{p}e^{2\int_0^{u}r(x)dx}du.\label{rho=inttegrale}
\end{equation}
$\rho$ is bijective if $r$ is integrable because $\rho$ is stricltly increasing and continous, hence $\rho$ define a bijection from $\mathbb{R}$ to $\rho(\mathbb{R})$ where 
$$\rho(\mathbb{R})=\left]\int_0^{-\infty}e^{2\int_0^{u}r(x)dx}du,\int_0^{+\infty}e^{2\int_0^{u}r(x)dx}du\right[$$
In conclusion for this type of function $\rho$, we have $\rho^{-1}$ is well defined, so we are able to transform the FBSDE \eqref{FBSDE in case of g=q+rZ^2} to a new FBSDE such that the generator of the BSDE does not depend on $z$

\begin{align}
    \begin{cases}
        dX_t &= b(t, X_t) dt + \sigma(t, X_t) dW_t, \\
        d\bar{Y}_t &= -\rho^{'}\circ\rho^{-1}(\bar{Y_t})q(t,X_t,\rho^{-1}(\bar{Y_t}))dt + Z_t dW_t, \\
        Y_T &= \rho\circ H(X_T).
    \end{cases}
\end{align}
Here, the Brownian motion, as well as $\sigma$ and $b$, remain unchanged. Thus, we return to the case where the generator $g$ does not depend on $z$, but with different terminal condition and different $Y$.
\section{Conclusion}
This study investigated the Lie symmetry of a new class of semi-linear PDE of the form \eqref{PDE-FBSDE}, which are connected to the FBSDE \eqref{FBSDE} through the generalized Feynman-Kac formula. Furthermore, leveraging the results of \cite{sym11091153}, we examined the symmetries of FBSDE \eqref{FBSDE}, investigating the interrelation between the symmetries of the FBSDE and its associated PDE. \\
Our analysis revealed conditions under which the PDE \eqref{PDE-FBSDE} admits Lie symmetries. Additionally, we showed that the set of symmetries of an FBSDE is contained within the set (algebra) of symmetries of the associated PDE, assuming that the generator $g$ of the FBSDE is independent of the variable $z$. Moreover, for certain specific (less general) cases of $g$ that depends on $z$, it is possible to transform the associated FBSDE  into a new FBSDE with a  generator that is independent of $z$, under some additional conditions.\\
It is also interesting to investigate the symmetries of FBSDE with jumps, as well as the symmetries the associated parabolic integro-differential PDE and the relation between these symmetries.\\
Consider a FBSDE with jump 
\begin{align}
    \begin{cases}
        X_t &= X_0+\displaystyle\int_0^t b(s, X_s) ds+\displaystyle\int_0^t \sigma(s, X_s) dB_s+\displaystyle\int_0^t\int_{\mathbb{R}}\gamma(s,X_{s^-},e)\Tilde{N}(ds,de) \\
        Y_{t} &=H(X_T)+\displaystyle\int_{t}^{T}g(s,X_s,Y_s,Z_s,U_s) ds-\displaystyle\int_{t}^{T}Z_s d B_{s}-\displaystyle\int_{t}^{T}\int_{\mathbb{R}}U_s(e)\Tilde{N}(ds,de),\quad t\in [0,T],\label{FBSDE with jumps}
    \end{cases}
\end{align}
where 
\begin{itemize}
    \item 
    $\gamma(s,X_{s^-},e)$ represents the size of the jump at time $s$ and state $X_{s^-}$(right before the jump).
    \item 
    $U_s(e)$ is the adapted process associated with the jump component.
\end{itemize}
The FBSDE \eqref{FBSDE with jumps} with jumps is connected to a parabolic integro-differential PDE through the Feynman Kac representation (Non linear case). The solution of the FBSDE $Y_t=u(t,X_t)$, where $u(t,X_t)$ is a deterministic function, satisfies the following parabolic PDE:
\begin{align}
    \begin{cases}
       &\displaystyle\frac{\partial u}{\partial t}(t,x)+\mathcal{L}u(t,x)+g(t,x,u,\sigma u_x,U(t,x))=0, \\
         &u(T,x)=g(x).\label{parabolic integro diff equa}
    \end{cases}
\end{align}
where:
\begin{itemize}
    \item 
    $\mathcal{L}$ is the infinitesimal generator of the forward process $X_t$, whcih includes both continous and jump dynamics. For a jump-diffusion, $\mathcal{L}$ is defined as:
    \begin{equation}
        \mathcal{L}u(t,x)=b(t,x)u_x(t,x)+\frac{1}{2}\sigma^2(t,x)u_{xx}(t,x)+\int_{\mathbb{R}}[u(t,x+\gamma(t,x,e))-u(t,x)-u_x(t,x)\gamma(t,x,e)]\nu(de)
    \end{equation}
\end{itemize}
When jumps are absent, i.e. $U_s=\gamma=0$, the FBSDE \eqref{FBSDE with jumps} reduces to the classical form without jump \eqref{FBSDE}, and the corresponding PDE to form \eqref{PDE-FBSDE}. One of the main contribution in this article is the study of symmetries of \eqref{FBSDE} and \eqref{PDE-FBSDE},  as well as the connection between the symmetries of the FBSDE and the associated PDE. A natural question that arises is the exploration of the symmetries of  \eqref{FBSDE with jumps} and \eqref{parabolic integro diff equa} and the investigation of the relationship between them.\\
\noindent  Corollary (8.5.5) in \cite[p.~157]{oksendal2013stochastic}  guarantee us that the transformed Brownian motion remains a Brownian motion. However, for jump processes like the Poisson process, do we have a similar theorem or a corollary analogous to the Brownian motion case ?\\
Aminu Nass, in \cite{nass2017w,nass2016n,nass2016symmetry}, attempted to investigate the invariance of the Poisson process, but only in a very specific case. \\
Another critical question is how the method of Lie symmetries of a PDE can be extended to integro-differential equations, specifically for an equation of the form \eqref{parabolic integro diff equa}. The primary challenge lies in handling the integral term within \eqref{parabolic integro diff equa}, as the classical approach presented in \cite{OLV86} cannot be directly applied. Nevertheless, Ibragimov, in \cite{ibragimov1995crc,grigoriev2010symmetries}, investigated the symmetries of certain integro-differential PDEs through alternative approaches. This raises the question: How can the symmetries of \eqref{parabolic integro diff equa} be effectively identified and characterized?
\bibliography{biblio2}

\begin{thebibliography}{BK13b}

\bibitem[Bah02]{10.1214/ECP.v7-1058}
Khaled Bahlali.
\newblock {Existence and Uniqueness of Solutions for BSDEs with Locally Lipschitz Coefficient}.
\newblock {\em Electronic Communications in Probability}, 7:169 -- 179, 2002.

\bibitem[BH08]{briand2008quadratic}
Philippe Briand and Ying Hu.
\newblock Quadratic bsdes with convex generators and unbounded terminal conditions.
\newblock {\em Probability Theory and Related Fields}, 141:543--567, 2008.

\bibitem[Bis73]{bismut1973conjugate}
Jean-Michel Bismut.
\newblock Conjugate convex functions in optimal stochastic control.
\newblock {\em Journal of Mathematical Analysis and Applications}, 44(2):384--404, 1973.

\bibitem[BK13a]{10.1214/12-AOP743}
Pauline Barrieu and Nicole~El Karoui.
\newblock {Monotone stability of quadratic semimartingales with applications to unbounded general quadratic BSDEs}.
\newblock {\em The Annals of Probability}, 41(3B):1831 -- 1863, 2013.

\bibitem[BK13b]{blu}
George~W Bluman and Sukeyuki Kumei.
\newblock {\em Symmetries and differential equations}, volume~81.
\newblock Springer Science \& Business Media, 2013.

\bibitem[Bow05]{Bowers}
Adam Bowers.
\newblock Classification of three-dimensional real lie algebras.
\newblock {\em available online at the webpage math. ucsd. edu/abowers/downloads/survey/3d\_Lie\_alg\_classify. pdf}, 2005.

\bibitem[Car16]{carmona2016lectures}
Ren{\'e} Carmona.
\newblock {\em Lectures on BSDEs, stochastic control, and stochastic differential games with financial applications}.
\newblock SIAM, 2016.

\bibitem[CD10]{craddock2010equivalence}
MJ~Craddock and Anthony~H Dooley.
\newblock On the equivalence of lie symmetries and group representations.
\newblock {\em Journal of Differential Equations}, 249(3):621--653, 2010.

\bibitem[CG20]{craddock2020lie}
Mark Craddock and Martino Grasselli.
\newblock Lie symmetry methods for local volatility models.
\newblock {\em Stochastic Processes and their Applications}, 130(6):3802--3841, 2020.

\bibitem[CL07]{craddock2007lie}
Mark Craddock and Kelly~A Lennox.
\newblock Lie group symmetries as integral transforms of fundamental solutions.
\newblock {\em Journal of Differential Equations}, 232(2):652--674, 2007.

\bibitem[CL09]{craddock2009calculation}
Mark Craddock and Kelly~A Lennox.
\newblock The calculation of expectations for classes of diffusion processes by lie symmetry methods.
\newblock 2009.

\bibitem[Gae00]{Giuseppe-Gaeta_2000}
Giuseppe Gaeta.
\newblock Lie-point symmetries and stochastic differential equations: Ii.
\newblock {\em Journal of Physics A: Mathematical and General}, 33(27):4883, jul 2000.

\bibitem[Gae19]{gaeta2019w}
Giuseppe Gaeta.
\newblock W-symmetries of ito stochastic differential equations.
\newblock {\em Journal of Mathematical Physics}, 60(5):053501, 2019.

\bibitem[GQ99]{Giuseppe-Gaeta_1999}
Giuseppe Gaeta and Niurka~Rodríguez Quintero.
\newblock Lie-point symmetries and stochastic differential equations.
\newblock {\em Journal of Physics A: Mathematical and General}, 32(48):8485, dec 1999.

\bibitem[Gri10]{grigoriev2010symmetries}
Yurii~N Grigoriev.
\newblock {\em Symmetries of Integro-Differential Equations}.
\newblock Springer, 2010.

\bibitem[GS17]{gaeta2017random}
Giuseppe Gaeta and Francesco Spadaro.
\newblock Random lie-point symmetries of stochastic differential equations.
\newblock {\em Journal of Mathematical Physics}, 58(5):053503, 2017.

\bibitem[Ibr95]{ibragimov1995crc}
Nail~H Ibragimov.
\newblock {\em CRC handbook of Lie group analysis of differential equations}, volume~3.
\newblock CRC press, 1995.

\bibitem[Jia09]{jia2009some}
Guangyan Jia.
\newblock Some uniqueness results for one-dimensional bsdes with uniformly continuous coefficients.
\newblock {\em Statistics \& probability letters}, 79(4):436--441, 2009.

\bibitem[Kob00]{10.1214/aop/1019160253}
Magdalena Kobylanski.
\newblock {Backward stochastic differential equations and partial differential equations with quadratic growth}.
\newblock {\em The Annals of Probability}, 28(2):558 -- 602, 2000.

\bibitem[Koz10]{kozlov2010symmetries}
Roman Kozlov.
\newblock Symmetries of systems of stochastic differential equations with diffusion matrices of full rank.
\newblock {\em Journal of Physics A: Mathematical and Theoretical}, 43(24):245201, 2010.

\bibitem[Koz20]{kozlov2020symmetries}
Roman Kozlov.
\newblock Symmetries of kolmogorov backward equation.
\newblock {\em arXiv preprint arXiv:2010.16128}, 2020.

\bibitem[Koz21]{kozlov2021symmetries}
Roman Kozlov.
\newblock Symmetries of kolmogorov backward equation.
\newblock {\em Journal of Nonlinear Mathematical Physics}, 28(2):182--193, 2021.

\bibitem[LQZ14]{lescot2014solving}
Paul Lescot, Helene Quintard, and Jean-Claude Zambrini.
\newblock Solving stochastic differential equations with cartan's exterior differential systems.
\newblock {\em arXiv preprint arXiv:1404.4802}, 2014.

\bibitem[LV19]{lescot2019symmetries}
Paul Lescot and Laur{\`e}ne Valade.
\newblock Symmetries of partial differential equations and stochastic processes in mathematical physics and in finance.
\newblock In {\em Journal of Physics: Conference Series}, volume 1194, page 012070. IOP Publishing, 2019.

\bibitem[NF16a]{nass2016n}
Aminu~M Nass and E~Fredericks.
\newblock N-symmetry of it{\^o} stochastic differential equation driven by poisson process.
\newblock {\em Int. J. Pure Appl. Math}, 110:165--182, 2016.

\bibitem[NF16b]{nass2016symmetry}
Aminu~M Nass and E~Fredericks.
\newblock Symmetry of jump-diffusion stochastic differential equations.
\newblock {\em Global and Stochastic Analysis}, 3(1):11--23, 2016.

\bibitem[NF17]{nass2017w}
Aminu~M Nass and E~Fredericks.
\newblock W-symmetries of jump-diffusion it{\^o} stochastic differential equations.
\newblock {\em Nonlinear Dynamics}, 90:2869--2877, 2017.

\bibitem[{\O}ks85]{oksendal1985stochastic}
Bernt {\O}ksendal.
\newblock When is a stochastic integral a time change of a diffusion?
\newblock {\em Preprint series: Pure mathematics http://urn. nb. no/URN: NBN: no-8076}, 1985.

\bibitem[Oks13]{oksendal2013stochastic}
Bernt Oksendal.
\newblock {\em Stochastic differential equations: an introduction with applications}.
\newblock Springer Science \& Business Media, 2013.

\bibitem[Olv86]{OLV86}
Peter~J. Olver.
\newblock {\em Applications of {Lie} groups to differential equations}, volume 107 of {\em Grad. Texts Math.}
\newblock Springer, Cham, 1986.

\bibitem[Ovs14]{ovsiannikov2014group}
Lev~Vasil'evich Ovsiannikov.
\newblock {\em Group analysis of differential equations}.
\newblock Academic press, 2014.

\bibitem[PP90]{PARDOUX199055}
E.~Pardoux and S.G. Peng.
\newblock Adapted solution of a backward stochastic differential equation.
\newblock {\em Systems $\&$ Control Letters}, 14(1):55--61, 1990.

\bibitem[SMS07]{srihirun2007definition}
Boonlert Srihirun, Sergey~V Meleshko, and Eckart Schulz.
\newblock On the definition of an admitted lie group for stochastic differential equations.
\newblock {\em Communications in Nonlinear Science and Numerical Simulation}, 12(8):1379--1389, 2007.

\bibitem[{\"U}na03]{unal2003symmetries}
Gazanfer {\"U}nal.
\newblock Symmetries of it{\^o} and stratonovich dynamical systems and their conserved quantities.
\newblock {\em Nonlinear Dynamics}, 32:417--426, 2003.

\bibitem[ZJ19]{sym11091153}
Na~Zhang and Guangyan Jia.
\newblock Lie-point symmetries and backward stochastic differential equations.
\newblock {\em Symmetry}, 11(9), 2019.

\bibitem[ZJ21]{ZHANG2021105527}
Na~Zhang and Guangyan Jia.
\newblock W-symmetries of backward stochastic differential equations, preservation of simple symmetries and kozlov’s theory.
\newblock {\em Communications in Nonlinear Science and Numerical Simulation}, 93:105527, 2021.

\end{thebibliography}

\bibliographystyle{alpha}

\end{document}